\documentclass[12pt]{article}
\usepackage{amsmath,amssymb}

\newfont{\zapfxii}{eusm10 scaled\magstep1}
\newcommand{\zap}[1]{\mbox{\zapfxii #1}}

\newcommand{\al}{\alpha}
\newcommand{\be}{\beta}
\newcommand{\ga}{\gamma}
\newcommand{\de}{\delta}
\newcommand{\ep}{\epsilon}
\newcommand{\la}{\lambda}
\newcommand{\ph}{\varphi}
\newcommand{\si}{\sigma}
\newcommand{\vep}{\varepsilon}

\newcommand{\De}{\Delta}
\newcommand{\Ga}{\Gamma}
\newcommand{\Si}{\Sigma}

\newcommand{\cN}{{\cal N}}
\newcommand{\cP}{{\cal P}}
\newcommand{\cR}{{\cal R}}
\newcommand{\cS}{{\cal S}}
\newcommand{\cT}{{\cal T}}

\def\fa{\mathfrak a}
\def\fg{\mathfrak g}
\def\fgl{\mathfrak{gl}}
\def\fosp{\mathfrak{osp}}
\def\fpl{\mathfrak{pl}}
\def\fs{\mathfrak{s}}
\def\fsl{\mathfrak{sl}}
\def\fspl{\mathfrak{spl}}
\def\fD{\mathfrak D}
\def\fP{\mathfrak P}
\def\fR{\mathfrak R}
\def\fS{\mathfrak S}
\def\fT{\mathfrak T}
\newcommand{\fU}{{\mathfrak U}}

\newcommand{\zS}{{\zap S}}
\newcommand{\zT}{{\zap T}}

\newcommand{\bcN}{\boldsymbol{\cN}}
\newcommand{\bcR}{\boldsymbol{\cR}}
\newcommand{\bcS}{\boldsymbol{\cS}}
\newcommand{\bcT}{\boldsymbol{\cT}}
\newcommand{\bfD}{\boldsymbol{\fD}}
\newcommand{\bfR}{\boldsymbol{\fR}}
\newcommand{\bfS}{\boldsymbol{\fS}}
\newcommand{\bfT}{\boldsymbol{\fT}}

\def\Je#1{J^\ep_{#1}}
\def\Jp#1{J^+_{#1}}
\def\J0#1{J^0_{#1}}
\def\Jm{J^-}
\def\Ke#1{K^\ep_{#1}}
\def\Kp#1{K^+_{#1}}
\def\K0#1{K^0_{#1}}
\def\Km{K^-}
\def\Qm{Q^-}
\def\Qp{Q^+}
\def\Qab{Q_{\al\be}}
\def\Se{S^\ep}
\def\Te{T^\ep}
\newcommand{\Tk}{T^{(k)}}

\def\bq{\,\overline{\!q}{}}
\def\brho{\bar\rho}
\def\tJ{\tilde J}
\def\tpi{\tilde\pi}

\def\tfD{\tilde\fD}
\def\tfR{\tilde\fR}
\def\tfS{\tilde\fS}
\def\tfT{\tilde\fT}
\def\tbfR{\boldsymbol{\tilde{\fR}}}
\def\tbfS{\boldsymbol{\tilde{\fS}}}
\def\tbfT{\boldsymbol{\tilde{\fT}}}

\newcommand{\CC}{{\mathbb C}}
\newcommand{\PP}{{\mathbb P}}
\newcommand{\RR}{{\mathbb R}}
\newcommand{\ZZ}{{\mathbb Z}}
\newcommand{\id}{1\hspace{-.25em}\rm{I}}

\newcommand{\Dx}{\partial _x}
\newcommand{\Dy}{\partial _y}

\newcommand{\iso}{\simeq}
\newcommand{\da}{\downarrow}
\newcommand{\ra}{\rightarrow}

\newcommand{\diag}{\operatorname{diag}}
\newcommand{\im}{\operatorname{im}}

\newcommand{\C}{\operatorname{C}}
\newcommand{\Deg}{\operatorname{Deg}}
\newcommand{\End}{\operatorname{End}}
\newcommand{\G}{\operatorname{G}}

\newcommand{\sch}{Schr\"o\-ding\-er}
\newcommand{\qes}{quasi-exactly solvable}
\newcommand{\Qes}{Quasi-exactly solvable}
\newcommand{\fdim}{fi\-nite-di\-men\-sion\-al}
\newcommand{\dfo}{differential operator}
\newcommand{\QED}{\quad Q.E.D.}

\def\cf{cf.~}

\newtheorem{thm}{\bf Theorem}[section]
\newtheorem{lemma}[thm]{\bf Lemma}

\newtheorem{cor}[thm]{\bf Corollary}

\numberwithin{equation}{section}

\newcommand{\bpm}{\begin{pmatrix}}
\newcommand{\epm}{\end{pmatrix}}
\def\ni{\noindent}
\def\ns{\normalsize}
\def\nem{\ns\em}

\newcounter{mylc}
\renewcommand{\themylc}{\roman{mylc}}
\newenvironment{mylist}{\begin{list}{\themylc )}
{\usecounter{mylc}\settowidth{\labelwidth}{iiii)}}}{\end{list}}

\title{\Large\bf{\mbox{\hspace{-1.2cm}The Lie Algebraic Structure
of Differential Operators}\\
Admitting Invariant Spaces of Polynomials\vspace{.5cm}}}

\author{\ns\textsc{Federico Finkel}\thanks{Supported
by DGICYT grant PB95--0401}\\
\nem Departamento de F\'\i sica Te\'orica II\\
\nem Universidad Complutense\\
\nem 28040 Madrid\\
\nem Spain
\and\ns\textsc{Niky Kamran}\thanks{On sabbatical leave from
the Department of Mathematics and Statistics, McGill University, Montr\'eal,
Qu\'ebec, Canada H3A 2K6. Supported by NSERC grant \#0GP0105490}\\
\nem The Fields Institute for Research\\
\nem in Mathematical Sciences\\
\nem Toronto, Ontario\\
\nem Canada M5T 3J1\vspace{.5cm}}

\date{\ns December 18, 1996}

\begin{document}
\maketitle

\begin{abstract}
We prove that the scalar and $2\times 2$
matrix \dfo s which preserve the simplest scalar and vector-valued
polynomial modules in two variables have a fundamental Lie algebraic
structure. Our approach is based on a general graphical method
which does not require the modules to be irreducible under
the action of the corresponding Lie (super)algebra. This method can be
generalized to modules of polynomials in an arbitrary number of variables.
We give generic examples of partially solvable \dfo s which are not Lie
algebraic. We show that certain vector-valued modules give rise to new
realizations of \fdim{} Lie superalgebras by first-order \dfo s.
\end{abstract}
\newpage

\section{Introduction}\label{sec.intro}

There has been a significant amount of interest recently in the
study of \sch{} operators which preserve explicit
\fdim{} subspaces $\cN$ of the underlying space of
smooth functions, \cite{Us94}, \cite{GKO94}, \cite{IL95}.
A \sch{} operator, or more generally any differential operator $T$,
which enjoys this property is said to be {\em partially solvable}. 
Indeed, if $\cN$ is $n$-dimensional and the restriction of $T$ to $\cN$ is a
self-adjoint operator with respect to a suitable inner product,
one can then determine $n$ eigenvalues and eigenfunctions
of $T$ algebraically, counting multiplicities. Notable examples of
partially solvable \sch{} operators in one dimension include some
new families of anharmonic oscillator potentials, as well as
the familiar Morse, P\"oschl-Teller and Mathieu potentials, which have been
thoroughly analyzed from this algebraic point of view, \cite{Tu88},
\cite{Sh89}, \cite{KO90}, \cite{GKO93}. It has been observed that
most of the known examples of partially solvable \sch{} operators can
be written as elements of the universal
enveloping algebra of a finite-dimensional
Lie algebra of first-order \dfo s admitting a
\fdim{} subspace of the underlying Hilbert space as a module of
smooth functions. The \dfo s satisfying this property are called {\em
\qes}, or more generally {\em Lie algebraic} if no
assumption is made on the existence of such a \fdim{} module.
One can thus construct many examples of
partially solvable \sch{} operators by first classifying
the Lie algebras of first-order \dfo s admitting
\fdim{} modules of functions, and then looking for elements in the
universal enveloping algebra which can be put in the \sch{} form
after conjugation by a suitable non-vanishing
multiplication operator and a local change of the independent variable.
This conjugation is referred to usually as a
gauge transformation. Note that the gauge transformation will
generally {\em not} be unitary, so that square integrability need
not hold at the outset, but only for the eigenfunctions of the
resulting \sch{} operator obtained after conjugation. This
program has proved to be remarkably successful in that it has enabled
one to construct many new examples of \qes{}
potentials in higher-dimensions, \cite{ST89}, \cite{GKO94}, as
well as matrix potentials for particles with spin, \cite{Sh89},
\cite{BK94}, \cite{FGR96}.

For all the examples known so far, the invariant
modules have a particularly simple structure when expressed in a
suitable coordinate system and gauge. In the one-dimensional case, they
are the vector spaces of polynomials of degree less or equal
than a positive integer $n$, and the corresponding Lie algebra is
simply the standard representation of $\fsl_2$ by first-order
\dfo s in one variable. It is easy to show that there are no
other possibilities  in one dimension, \cite{Mi68}, \cite{KO90}.
In the case of \qes{} \dfo s in two variables associated to a complex
Lie algebra, the modules also take the form of polynomial modules
in the appropriate coordinates and gauge, \cite{GKO91}.
They are either one of the vector spaces of bivariate polynomials
$x^iy^j$ of bidegree $(i,j)$ less or equal than an ordered pair $(n,m)$
of positive integers, or one of the vector spaces of polynomials
of total degree $i+j$ less or equal than a positive
integer $n$, or one of the vector spaces of polynomials whose
bidegree $(i,j)$ satisfies a linear constraint of the form $i+rj \leq p$,
for some positive integers $r$ and $p$. We will refer to these
modules as the rectangular, triangular or staircase modules in view of
their structure within the lattice of bidegrees. The corresponding Lie algebras
are given by the standard realizations of $\fsl_2\oplus\fsl_2$,  $\fsl_3$ and
$\fgl_2\ltimes\CC^{r+1}$ by first-order \dfo s in the plane. The
inequivalent real forms of these algebras give rise to different
quasi-exactly solvable potentials, \cite{GKO96}. 
It should be noted that while the rectangular and triangular modules
are acted on irreducibly by $\fsl_2\oplus\fsl_2$ and $\fsl_3$,
this is {\em not} the case for the staircase module.

A fundamental question is to determine why practically every known example
of partially solvable \dfo{} should be expressible as a polynomial in
the generators of a \fdim{} Lie algebra of first-order \dfo s leaving
the corresponding subspace invariant, and to what extent such a polynomial
is unique. In particular, one would like to know if there exist
modules for which the generic partially solvable \dfo s are not Lie algebraic.
These are the main problems that we address in this paper for the case
of rectangular, triangular and staircase scalar or rank 2 modules.
In the case in which the module is acted on irreducibly by the Lie
algebra, which in our case would be the rectangular and triangular modules,
a partial answer is provided by the Burnside Theorem, as pointed out by
Turbiner, \cite{Tu94}. Recall that according to the Burnside Theorem, if $V$
is a complex vector space and $\fa$ is a subalgebra of $\End V$ acting
irreducibly, then every endomorphism of
$V$ can be expressed as a polynomial in the generators of $\fa$.
This implies that any \dfo{} preserving a rectangular or a triangular
module can be expressed as the sum of a polynomial in the generators
of the realization of the corresponding Lie algebra
and a \dfo{} annihilating the module. Note, however, that this argument
does not provide any direct
information on the relation between the order of a \dfo{}
leaving a rectangular or triangular module invariant and the minimum
degree of a polynomial representation thereof,
owing to the presence of polynomial
relations between the generators of the Lie algebra, which form a
basis for the primitive ideal associated to the
module. The case of the staircase modules is considerably more
complicated, since the underlying Lie algebras of first-order
\dfo s do not act irreducibly, so that one cannot apply the Burnside Theorem.
In fact, we will see that in stark contrast with the rectangular and 
triangular cases, the generic differential operator preserving the staircase
module is {\em not} Lie-algebraic. One needs therefore to 
develop a different approach altogether. This is what
we suceed to do in this paper,
by developing a general graphical method which does not require the
(scalar or rank 2) modules to be irreducible under
the action of the corresponding algebra. We will 
thus entirely bypass the use of the Burnside Theorem and the need for
extensive calculations in local coordinates. Indeed, the latter can get
prohibitively difficult in higher dimensions. (Even
the one-dimensional scalar case is quite involved when one works in
coordinates, \cite{Tu92}). It is worth noting that the graphical method can
be generalized without  major difficulties to the case of rank $k$ modules of
polynomials in $N$ variables. It is thus a powerful tool for
future applications.

In Section~2, we briefly review the normal forms for the
maximal Lie algebras of first-order \dfo s in two
complex variables which admit finite-dimensional modules of
polynomials when expressed in a suitable coordinate
system and gauge. In Section~3, we
unravel by means of a general graphical method the Lie algebraic
structure of the scalar \dfo s of arbitrary finite order which preserve a
rectangular, triangular or staircase module.
We give a formula for the minimum degree of a polynomial in the generators
of the corresponding Lie algebra representing the given \qes{} \dfo.
In the case of triangular modules, we also obtain an exact
formula for the number of free
parameters which determine the most general partially solvable \dfo{}
of arbitrary finite order. Section~4 is concerned with 
$2\times 2$ matrix \dfo s of arbitrary finite order admitting the
direct sum of two rectangular, triangular or staircase modules as an
invariant subspace. For suitable values of the gap parameters, the underlying
algebras of matrix \dfo s give rise to new realizations of \fdim{} Lie
superalgebras of first-order \dfo s. The abstract structure
of these \fdim{} Lie superalgebras is identified in every case.
For the remaining values of the gap
parameters, one obtains an infinite-dimensional Lie superalgebra.
Remarkably, the graphical method introduced for scalar modules
enables us in both the finite and infinite-dimensional cases
to prove structure theorems for the corresponding partially solvable
matrix \dfo s.

\section{\Qes{} Lie algebras in \boldmath$\CC^2$}\label{sec.qes}

Our purpose in this section is to briefly recall from \cite{GKO91}
the normal forms for the maximal Lie algebras of first-order \dfo s
in two complex variables admitting
\fdim{} modules of polynomials in suitable coordinates,
along with the structure of these modules.

Let $U$ be an open subset of $\CC^2$, with local coordinates $(x,y)$.
The first-order \dfo s on $U$,
$$
T=f(x,y)\Dx+g(x,y)\Dy+h(x,y)\,,
$$
with coefficients in $\C^\infty(U,\CC)$ form a Lie algebra which we
will denote by $\fD^{(1)}=\fD^{(1)}(U)$. The space $\C^\infty(U,\CC)$ is
naturally a $\fD^{(1)}$-module. We will use the terminology of \cite{GKO94},
and refer to \fdim{} subalgebras $\fg$ of $\fD^{(1)}$
admitting a \fdim{} module $\cN\subset\C^\infty(U,\CC)$ as {\em\qes{}}
Lie algebras. The local diffeomorphisms $\ph$ of $U$ and the rescalings
by smooth functions $e^\si\in\C^\infty(U,\CC^\ast)$ define an
infinite pseudogroup $\G$ of automorphisms
$\Phi$ of $\fD^{(1)}$, whose action is given by
\begin{equation}\label{Phi}
\Phi(T)=e^\si\cdot\ph_\ast T\cdot e^{-\si}\,.
\end{equation}
The classification under the action of $\G$ of all \qes{} Lie algebras
in two complex variables is known, \cite{GKO91} (its real counterpart
has been recently completed, \cite{GKO96}). We shall focus our attention
on the three maximal families of equivalence classes of \qes{} Lie
algebras whose normal forms preserve a module of polynomials.

\begin{mylist}
\item $\fg^{11}_{n,m}\iso\fsl_2\oplus\fsl_2$, with generators:
\begin{gather}
\hspace{-.9cm}\Jp n=x^2\Dx-n x\,,\qquad \Jm=\Dx\,,\qquad
\J0 n=x\Dx-\frac n2\,,\notag\\
\Kp m=y^2\Dy-m y\,,\qquad \Km=\Dy\,,\qquad \K0 m=y\Dy-\frac m2\,,\label{g11}
\end{gather}
with $n,m\in\ZZ^+$, and associated \fdim{} module given by
$$
\cR_{n,m}=\{x^iy^j\:|\:0\leq i\leq n\,,\; 0\leq j\leq m\}\,.
$$
\item $\fg^{15}_{n}\iso\fsl_3$, with generators:
\begin{gather}
J^1_n=x^2\Dx+xy\Dy-nx\,,\qquad J^2_n=xy\Dx+y^2\Dy-ny\,,\notag\\[.3cm]
J^3=y\Dx\,,\qquad J^4=\Dx\,,\qquad J^5=\Dy\,,\qquad J^6=x\Dy\,,\notag\\[.1cm]
J^7_n=x\Dx-\frac n3\,,\qquad J^8_n=y\Dy-\frac n3\,,\label{g15}
\end{gather}
with $n\in\ZZ^+$ and corresponding \fdim{} module given by
$$
\cT_n=\{x^iy^j\:|\:0\leq i+j\leq n\}\,.
$$
\item $\fg^{24,r}_{p}\iso\fgl_2\ltimes\CC^{r+1}$, with generators:
\begin{gather}
J^1=\Dx\,,\qquad J^2_p=x^2\Dx+r x y\Dy-p x\,,\qquad
J^3_p=x\Dx-\frac p2\,,\notag\\[.1cm]
J^4=y\Dy\,,\qquad J^{5+i}=x^i\Dy\,,\qquad 0\leq i\leq r\,,\label{g24}
\end{gather}
with $2\leq r\in\ZZ^+$, $p\in\ZZ^+$ and, for $q\in\ZZ^+$,
associated module given by
$$
\cS^r_{p,q}=\{x^iy^j\:|\:0\leq i+rj\leq p\,,\; 0\leq j\leq q\}\,,
$$
As a special case, we shall consider the module $\cS^r_p=\cS^r_{p,[p/r]}$.
\end{mylist}
The nomenclature of the algebras is the same as the one used in \cite{GKO91},
and the letters $\cR$, $\cT$ and $\cS$ used to denote the modules reflect
of course the rectangular, triangular and staircase structure of the
lattices in the positive $\ZZ^+\times\ZZ^+$ quadrant corresponding to the
powers $(i,j)$ that appear in $\cR_{n,m}$, $\cT_n$ and $\cS^r_{p,q}$. The
integers $n$, $m$, $p$, $q$ which label the realizations of
$\fsl_2\oplus\fsl_2$, $\fsl_3$ and $\fgl_2\ltimes\CC^{r+1}$ in the list
above are related to the global models for these algebras, where the
underlying complex surfaces are respectively given by
$\PP_1(\CC)\times\PP_1(\CC)$, $\PP_2(\CC)$ and the
$r$-th Hirzebruch surface $\Si_r$, and the modules are sections of holomorphic
line bundles over these surfaces, \cite{GHKO93}.

\section{Scalar $\mathbf k$-th order \dfo s in \boldmath$\CC^2$ with
invariant modules of polynomials}\label{sec.scalar}

Our aim in this section is to prove structure theorems for the subalgebras
of the associative algebra $\fD=\fD(U)$ of scalar \dfo s in two
complex variables
\begin{equation}\label{T}
T=\sum_{0\leq i,j<\infty}f_{ij}(x,y)\,\Dx^i\Dy^j\,,\qquad
f_{ij}(x,y)\in\C^\infty(U,\CC)\,,
\end{equation}
which map the modules $\cR_{n,m}$, $\cT_n$ and $\cS^r_{p,q}$ to
themselves. Let us remark that our results will also be valid when
$x$ and $y$ are restricted to
the real plane $\RR^2$. The motivation for our study arises from
a theorem of Turbiner concerning scalar $k$-th order
\dfo s in one variable, \cite{Tu92}. To set Turbiner's theorem in context,
let us recall, \cite{Mi68}, \cite{KO90}, that up to equivalence under the
one-dimensional analog to the pseudogroup $\G$ defined in~\eqref{Phi},
there is a unique family of \qes{} Lie algebras of first-order \dfo s on
the line, namely the family of representations of $\fsl_2$ given by
$\fg_n=\{\Jp n,\,\Jm,\,\J0 n\}$, where
\begin{equation}\label{Js}
\Jp n=x^2\Dx-nx\,,\qquad\Jm=\Dx\,,\qquad\J0 n=x\Dx-\frac n2\,.
\end{equation}
The associated invariant module is the polynomial module
$$
\cP_n=\{x^i\:|\:0\leq i\leq n\}\,.
$$
Let $\tfD$ be the associative algebra of all linear \dfo s in one complex
variable, and let $\tfD^{(k)}$ denote its subspace of \dfo s of order at
most $k$. We shall denote by
$\tfD^{(>k)}$ the subalgebra of $\tfD$ spanned by all \dfo s of the form
$T\,\Dx^{k+1}$, where $T$ is an arbitrary \dfo. Let $\fP_n$ be the
subalgebra of $\tfD$ of \dfo s which preserve $\cP_n$, and let
$\fP^{(k)}_n=\fP_n\cap\tfD^{(k)}$. Finally, let
$\zS^k(\fsl_2)=\bigoplus_{i=0}^k S^i(\fsl_2)$  be the direct sum of the
first $k$ symmetric powers of $\fsl_2$. Turbiner's theorem may be
reformulated as follows:
\begin{thm}\label{thm.Tu}
i) The linear map $\rho_n^k:\zS^k(\fsl_2)\ra\fP^{(k)}_n$ determined by the
representation $\rho_n:\fsl_2\ra\fg_n$ is surjective
for all $k=0,\ldots,n$.\quad
ii) If $k>n$, then $\fP^{(k)}_n=\fP^{(n)}_n\oplus(\tfD^{(k)}\cap\tfD^{(>n)})$.
\end{thm}
{\em Remark.} Turbiner's original proof was done by explicit computation.
We present here a simplified proof whose principle will serve as a guide
in the higher dimensional case.

{\ni\em Proof.} We first note that $\ker\pi_n=\tfD^{(>n)}$, where
$\pi_n$ denotes the induced homomorphism $\pi_n:\fP_n\ra\End(\cP_n)$,
so {\em ii)} clearly holds.
Therefore, the restriction of $\pi_n$ to 
$\fP^{(n)}_n$ is injective, so $\dim\fP^{(n)}_n\leq(n+1)^2$.
Now, the monomials $\{(J^\pm_n)^i(\J0 n)^{j-i}\}_{i=0}^j$, $j=0,1,\ldots$,
are all linearly independent \dfo s. In particular, there are
$(n+1)^2$ monomials of order at most $n$, so they form a basis of
$\fP^{(n)}_n$. Since $(J^\pm_n)^i(\J0 n)^{j-i}=x^{j\pm i}\Dx^j+\cdots$,
no monomial of degree higher than $k$ may appear in the expansion of a
\dfo{} $\Tk\in\fP^{(k)}_n$.\QED\\

Thus, Theorem~\ref{thm.Tu} simply expresses the fact that any \dfo
$$
T=\sum_{i=0}^{k\leq n} f_i(x)\,\Dx^i\,,
$$
preserving $\cP_n$ may be written as a polynomial of degree $k$ in
the generators~\eqref{Js}.

Later, Turbiner proposed a proof of Theorem~\ref{thm.Tu} based on the
irreducibility of the module $\cP_n$ under the action of $\fg_n$ and the
Burnside Theorem. However, this argument only guarantees that $\Tk$ will
be the sum of a polynomial in the generators of $\fg_n$ with an operator
annihilating the module. In particular, it gives no information on the
degree of the polynomial.\\

The linear map $\rho_n^k:\zS^k(\fsl_2)\ra\fP^{(k)}_n\subset\tfD^{(k)}$
is clearly not injective for \mbox{$k\geq 2$}. In fact, the
set $\{(J^\pm_n)^i(\J0 n)^{j-i}\}_{i=0}^j$, $j=0,\ldots,k$
is a basis of $\im\rho_n^k$ (the completeness follows from the quadratic
relation $\Jp n\Jm=\J0 n\J0 n-\J0 n-\frac 14n(n+2)$ which expresses the scalar 
action of the Casimir operator on $\cP_n$). Thus,
\begin{equation}\label{kerrho}
\hspace{-.2cm}\dim(\ker\rho_n^k)=\dim\zS^k(\fsl_2)-(k+1)^2=
\frac 16(k+1)k(k-1)\,.
\end{equation}
The dimension of the kernel of the induced map
$\brho_n^k:S^k(\fsl_2)\ra\tfD^{(k)}/\tfD^{(k-1)}$ may also be easily obtained:
\begin{equation}\label{kerbrho}
\dim(\ker\brho_n^k)=\frac 12 k(k-1)\,.
\end{equation}
We emphasize that equations \eqref{kerrho} and \eqref{kerbrho} are
still valid if the cohomology parameter $n$ labeling the representation
$\fg_n$ is replaced by an arbitrary complex number $\la$;
see~\cite{PV96} and~\cite{Di73} for more details.\\

We now focus on the two-variable modules $\cR_{n,m}$, $\cT_n$ and
$\cS^r_{p,q}$, and define $\fR_{n,m}$, $\fT_n$ and $\fS^r_{p,q}$
as the subalgebras of $\fD$ which preserve these
modules. Let $\pi_{n,m}:\fR_{n,m}\ra\End(\cR_{n,m})$,
$\pi_n:\fT_n\ra\End(\cT_n)$, and
\mbox{$\pi^r_{p,q}:\fS^r_{p,q}\ra\End(\cS^r_{p,q})$} denote the
induced homomorphisms of associative algebras. We have the following
elementary result:
\begin{lemma}\label{lemma.ker}
The kernels of $\pi_{n,m}$, $\pi_n$, and $\pi^r_{p,q}$ are given by
\begin{align*}
& \ker\pi_{n,m}=\big\{T\in\fD\;|\; T=\sum_{\substack{i>n\text{ \rm or}\\ j>m}}
   f_{ij}(x,y)\,\Dx^i\Dy^j\big\}\,,\\
& \hspace{2cm}\ker\pi_{n}=\big\{T\in\fD\;|\; T=\sum_{i+j>n}
   f_{ij}(x,y)\,\Dx^i\Dy^j\big\}\,,\\[.2cm]
& \hspace{4cm}\ker\pi^r_{p,q}=\big\{T\in\fD\;|\;
   T=\sum_{\substack{i+rj>p\\ \text{\rm or }j>q}}
   f_{ij}(x,y)\,\Dx^i\Dy^j\big\}\,.
\end{align*}
\end{lemma}
We now choose the following distinguished complements to the kernels
in Lemma~\ref{lemma.ker}:
\begin{align*}
& \tfR_{n,m}=\big\{T\in\fR_{n,m}\;|\;
     T=\sum_{\substack{0\leq i\leq n\\ 0\leq j\leq m}}
     f_{ij}(x,y)\,\Dx^i\Dy^j\big\}\,,\\
& \hspace{2cm}\tfT_n=\big\{T\in\fT_n\;|\; T=\sum_{0\leq i+j\leq n}
     f_{ij}(x,y)\,\Dx^i\Dy^j\big\}\,,\\[.2cm]
& \hspace{4cm}\tfS^r_{p,q}=\big\{T\in\fS^r_{p,q}\;|\;
     T=\sum_{\substack{0\leq i+rj\leq p\\0\leq j\leq q}}
     f_{ij}(x,y)\,\Dx^i\Dy^j\big\}\,.
\end{align*}
Note that the order of a \dfo{} in the distinguished complement
$\tfS^r_{p,q}$ cannot exceed $p$. The next result follows from a
simple constructive argument:
\begin{lemma}\label{lemma.iso}
The restrictions $\tpi_{n,m}:\tfR_{n,m}\ra\End(\cR_{n,m})$,
$\tpi_{n}:\tfT_n\ra\End(\cT_n)$,
and \mbox{$\tpi^r_{p,q}:\tfS^r_{p,q}\ra\End(\cS^r_{p,q})$}
are vector space isomorphisms.
\end{lemma}
{\em Proof.} Let us prove the lemma for one of the modules,
say $\cT_n$. Let $\zT\in\End(\cT_n)$. The \dfo{} $T\in\tfT_n$ with
coefficients $f_{ij}$ given by
$$
f_{00}=\zT(1)\,,\qquad f_{10}=\zT(x)-xf_{00}\,,\qquad
f_{01}=\zT(y)-yf_{00}\,,\ldots
$$
clearly satisfies $\tpi_{n}(T)=\zT$.\QED\\

Note that the coefficients $f_{ij}$ of the \dfo s in the
distinguished complements are polynomials in $x$ and $y$.
The subalgebra $\fD_P$ of $\fD$ of \dfo s~\eqref{T} with
polynomial coefficients $f_{ij}$ inherits a natural $\ZZ\times\ZZ$ grading
from $\CC[x,y]$. We shall say that a \dfo{} $T\in\fD_P$ has {\em bidegree}
$\deg T=(i,j)$ if $T(x^{i_0}y^{j_0})\in\langle x^{i_0+i}y^{j_0+j}\rangle$
for all $(i_0,j_0)\in\ZZ^+\times\ZZ^+$.
Note that the generators of $\fg^{11}_{n,m}$, $\fg^{15}_n$ and
$\fg^{24,r}_p$ all have well-defined
bidegree, \cite{Tu92b}. A set of \dfo s in $\fD_P$ with different
bidegrees is of course linearly independent.

We will obtain in what follows structure theorems for each
of the distinguished complements. Let us remark that both $\fg^{11}_{n,m}$
and $\fg^{15}_n$ act irreducibly on their associated modules $\cR_{n,m}$
and $\cT_n$. Therefore, a direct application of the Burnside Theorem shows
that any \dfo{} $T$ in $\tfR_{n,m}$ (respectively $\tfT_n$) may be
constructed as a polynomial in the generators of $\fg^{11}_{n,m}$
(respectively $\fg^{15}_n)$ plus an element of $\ker\pi_{n,m}$
(respectively $\ker\pi_n$), \cite{Tu94}. Note that the action of
$\fg^{24,r}_p$ on $\cS^r_{p,q}$ is {\em neither} irreducible {\em nor}
completely reducible, contrary to the the assertion in ref.~\cite{Tu92b}.
In fact, we will exhibit \dfo s in
$\tfS^r_{p,q}$ which cannot be expressed as a polynomial in the
generators of $\fg^{24,r}_p$ plus an element in $\ker\pi^r_{p,q}$.

\subsection{The rectangular module $\cR_{n,m}$}

Let $\fD^{(k,l)}$ denote the subspace of $\fD$ spanned by all \dfo s
of $x$-order at most $k$ and $y$-order at most $l$, and let
$\tfR_{n,m}^{(k,l)}=\tfR_{n,m}\cap\fD^{(k,l)}$. The
structure theorem for $\tilde\fR_{n,m}$ is a
straightforward generalization of Theorem~\ref{thm.Tu}:
\begin{thm}\label{thm.R}
Let $\brho^i_n(x):S^i(\fsl_2)\ra\fP_n(x)$ and
$\brho^j_m(y):S^j(\fsl_2)\ra\fP_m(y)$ be the linear
maps determined by $\rho_n(x):\fsl_2\ra\fg_n(x)$
and $\rho_m(y):\fsl_2\ra\fg_m(y)$, respectively. Then,
$$
\bigoplus_{\substack{0\leq i\leq k\\0\leq j\leq l}}
\brho^i_n(x)\otimes\brho^j_m(y):
\bigoplus_{\substack{0\leq i\leq k\\0\leq j\leq l}}
S^i(\fsl_2)\otimes S^j(\fsl_2)\ra\tfR^{(k,l)}_{n,m}\,,
$$
is surjective for all $k=0,\ldots,n$ and $l=0,\ldots,m$.
\end{thm}

\subsection{The triangular module $\cT_n$}

Let $\fD^{(k)}$ be the subspace of $\fD$ spanned by all \dfo s of
order at most $k$, and let
$\tfT^{(k)}_n=\tfT_n\cap\fD^{(k)}$. As before, let
$\zS^k(\fsl_3)=\bigoplus_{i=0}^k S^i(\fsl_3)$
be the direct sum of the first $k$ symmetric powers of $\fsl_3$, and let
$\rho_n^k:\zS^k(\fsl_3)\ra\fD^{(k)}$ be the linear map determined by the
representation $\rho_n:\fsl_3\ra\fg^{15}_n$. Our purpose is to construct
a suitable basis of $\im\rho^k_n$ for all $k\in\ZZ^+$, from which the
structure theorem for $\tfT_n$ will arise as a simple corollary. The
bidegree of the generators \eqref{g15} of $\fg^{15}_n$ is given by:
\begin{gather*}
\deg J^1_n=(1,0)\,,\quad \deg J^2_n=(0,1)\,,\quad
\deg J^3=(-1,1)\,,\quad \deg J^4=(-1,0)\,,\\
\deg J^5=(0,-1)\,,\quad \deg J^6=(1,-1)\,,\quad
\deg J^7_n=(0,0)\,,\quad \deg J^8_n=(0,0)\,,
\end{gather*}
and the bidegree of a monomial $J^K_n=(J^1_n)^{k_1}\cdots(J^8_n)^{k_8}$
is simply $\sum_{i=0}^6 k_i\deg J^i_n$. Let $|K|=k_1+\cdots+k_8$ denote
the {\em degree} of $J^K_n$. Note that $|K|$ is just the order of $J^K_n$.
We define the {\em length} of a monomial $J^K_n$ as
$$
|J^K_n|=\min_{\deg J^L_n=\deg J^K_n}|L|\,,
$$
where the minimum is taken over the set of monomials of degree at most $|K|$.
In other words, if $\deg J^K_n=(i,j)$, then $|J^K_n|$ is just the minimum
number of ``steps'' required to map a point $(i_0,j_0)\in\ZZ^+\times\ZZ^+$
to \mbox{$(i_0+i,j_0+j)$}. A monomial $J^K_n$ is of
{\em maximal length} if $|J^K_n|=|K|$. The monomials $1,J^1_n,\ldots,J^6_n$
are obviously of maximal length. The following facts follow from an easy 
graphical argument, \cf Fig.~1:
\begin{mylist}
\item There is exactly one monomial of maximal length of bidegree $(i,j)$
for each $(i,j)\in\ZZ\times\ZZ$, up to the ordering of its factors.
\item The monomials of maximal length $l$ are
$\{(J^{s-1}_n)^j\,(J^s_n)^{l-j}\}_{j=1-\de_{l0}}^l$, with $s=1,\ldots,6$
and $J^0_n=J^6_n$.
\end{mylist}
\begin{figure}[h]
\begin{center}
\begin{picture}(240,80)(-120,-40)
\put(0,0){\circle*{3}}\thicklines\put(0,0){\circle{6}}\thinlines
\put(20,0){\circle{4}}\put(0,20){\circle{4}}\put(-20,20){\circle{4}}
\put(-20,0){\circle{4}}\put(0,-20){\circle{4}}\put(20,-20){\circle{4}}
\put(40,0){\circle*{4}}\put(20,20){\circle*{4}}\put(0,40){\circle*{4}}
\put(-20,40){\circle*{4}}\put(-40,40){\circle*{4}}\put(-40,20){\circle*{4}}
\put(-40,0){\circle*{4}}\put(-20,-20){\circle*{4}}\put(0,-40){\circle*{4}}
\put(20,-40){\circle*{4}}\put(40,-40){\circle*{4}}\put(40,-20){\circle*{4}}
\put(120,20){\circle*{3}}\thicklines\put(120,20){\circle{6}}\thinlines
\put(130,17){\small Length $l=0$}
\put(120,0){\circle{4}}\put(130,-3){\small Length $l=1$}
\put(120,-20){\circle*{4}}\put(130,-23){\small Length $l=2$}
\put(-11,-10){\footnotesize $(0,0)$}
\end{picture}
\end{center}
\caption{Lattice of bidegrees of monomials of length $l\leq 2$.}
\end{figure}
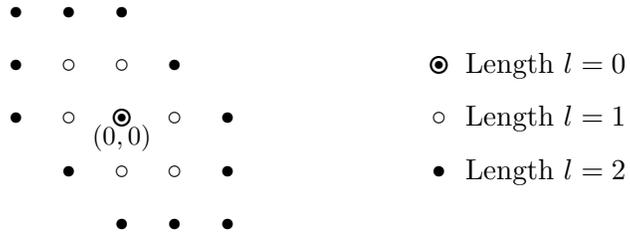
Note that there are $6l$ monomials of maximal length $l\geq 1$.
If $J^L_n$ is a monomial of maximal length,
the monomials $J^L_n\,(J^7_n)^i\,(J^8_n)^j$, $i,j\geq 0$,
are linearly independent, and have the same bidegree as $J^L_n$.
Remarkably, the monomials of this form are also a generating set of
$\im\rho^k_n$ for all $k\in\ZZ^+$, according to the following lemma:
\begin{lemma}\label{lemma.basis}\ \\[.3cm]
i) The monomials
\begin{equation}\label{basis}
\big\{(J^{s-1}_n)^j(J^s_n)^{l-j}(J^7_n)^i(J^8_n)^{t-l-i}\:\big|\:
1-\de_{l0}\leq j\leq l\leq t\leq k,\; 0\leq i\leq t-l\big\},
\end{equation}
with $s=1,\ldots,6$ and $J^0_n=J^6_n$,
form a basis of $\im\rho^k_n$ for all $k$.\\

\ni ii) The number of $k$-th order monomials in the basis is $n_k=(k+1)^3$,
and
\begin{equation}\label{dimim}
\dim(\im\rho^k_n)=\frac 14(k+1)^2(k+2)^2.
\end{equation}
\end{lemma}
{\ni\em Proof.} Linear independence of the monomials follows from the
preceeding remarks. Completeness is a consequence of the following 9
quadratic relations,
\begin{align}\label{rel}
& J^1_n J^3=J^2_n\big(J^7_n+\frac n3\big)\,,\qquad J^1_n J^4=
(J^7_n)^2+J^7_n J^8_n-J^7_n+\frac n3 J^8_n-\frac n3(\frac n3+1)\,,\notag\\
& J^1_n J^5=J^6\big(J^7_n+J^8_n-(\frac n3+1)\big)\,,\qquad
J^2_n J^4=J^3\big(J^7_n+J^8_n-(\frac n3+1)\big)\,,\notag\\
& J^2_n J^5=
J^7_n J^8_n+(J^8_n)^2+\frac n3 J^7_n-J^8_n-\frac n3(\frac n3+1)\,,\qquad
J^2_n J^6=J^1_n\big(J^8_n+\frac n3\big)\,,\notag\\
& J^3 J^5=J^4\big(J^8_n+\frac n3\big)\,,\qquad J^3 J^6=
J^7_n J^8_n+\frac n3 J^7_n+(\frac n3+1) J^8_n+\frac n3(\frac n3+1)\,,\notag\\
& J^4 J^6=J^5\big(J^7_n+\frac n3+1\big)\,,\
\end{align}
which allow us to reduce any monomial $J^K_n$ to a linear combination of
the elements of the set~\eqref{basis}. The second assertion now follows
by straightforward computation.\QED\\

From Lemma~\ref{lemma.basis}, we obtain the following result:
\begin{thm}\label{thm.T}
The linear map $\rho_n^k:\zS^k(\fsl_3)\ra\tfT^{(k)}_n$ determined by the
representation $\rho_n:\fsl_3\ra\fg^{15}_n$ is surjective for all
$k=0,\ldots,n$.
\end{thm}
{\ni\em Proof.} The map $\rho_n^n:\zS^n(\fsl_3)\ra\tfT_n$ is clearly
surjective, since $\dim\tfT_n=\frac 14(n+1)^2(n+2)^2$. The highest
order terms of the monomials in~\eqref{basis} are all independent,
which implies that no monomial of degree higher
than $k$ may appear in the expansion of a \dfo{} $\Tk\in\tfT^{(k)}_n$.\QED\\

Theorem~\ref{thm.T} is stated in a different form in~\cite{Tu92b}. It
should be noted, however, that just as in the $\fsl_2$ case, this
theorem does not follow directly
from the Burnside Theorem, contrary to what is stated in~\cite{Tu92b}.
As another immediate consequence of Lemma~\ref{lemma.basis}, we have:
\begin{cor}
The kernel of the homomorphism $\brho_n:\fU(\fsl_3)\ra\fD$ determined by the
representation $\rho_n:\fsl_3\ra\fg^{15}_n$ is the primitive ideal generated
by the quadratic relations~\eqref{rel}.
\end{cor}

Let us emphasize that these results are completely independent
of the discrete character of the cohomology parameter $n$ labeling the
Lie algebras $\fg^{15}_n$, being also valid for a complex
cohomology parameter $\la$, as considered
in \cite{Di73} and \cite{PV96}. The graphical method used to obtain
the basis~\eqref{basis} can be generalized to the \mbox{$N$-dimensional}
representations of $\fsl_{N+1}$ which
appear as the hidden symmetry algebra of the Calogero model and its \qes{}
extensions, \cite{Ca69}, \cite{RT95}, \cite{MRT96}.
Finally, let us give the formula analogous to~\eqref{kerbrho} for the map
$\brho^k_n:S^k(\fsl_3)\ra\fD^{(k)}/\fD^{(k-1)}$:
$$
\dim(\ker\brho_n^k)=\frac{(k+1)k(k-1)}{7!}
\big(k^4+28 k^3+323 k^2+1988 k+2052\big)\,.
$$

\subsection{The staircase module $\cS^r_{p,q}$}

The structure of the space of operators preserving this type of modules
is considerably more complicated than that of the other two modules.
As remarked before, the action of $\fg^{24,r}_p$ on $\cS^r_{p,q}$ is
reducible, leading to a number of interesting consequences.

The bidegrees of the generators of $\fg^{24,r}_p$ are:
\begin{align*}
& \deg J^1=(-1,0)\,,\qquad\deg J^2_p=(1,0)\,,\qquad\deg J^3_p=(0,0)\,,\\
& \deg J^4=(0,0)\,,\qquad \deg J^{5+i}=(i,-1)\,,\qquad 0\leq i\leq r\,.
\end{align*}
We shall refer to the second component of the bidegree as the $y$-degree.
Note that the $y$-degree of the generators~\eqref{g24} of $\fg^{24,r}_p$
is nonpositive. Thus, no operator in the distinguished complement
$\tfS^r_{p,q}$ with positive $y$-degree (as $T=y\Dx^2$,
if $r=2$, $p\geq 2$ and $q=[p/2]$), can be obtained as a polynomial
in the generators plus an operator in $\ker\pi^r_{p,q}$. Moreover,
not every operator in $\tfS^r_{p,q}$
with nonpositive $y$-degree may be written as a polynomial in the
generators of $\fg^{24,r}_p$ only. An example
is given by $T=x^2(x^2\Dx^2-2x\Dx-2y\Dy+2)$ in the case $p=2q=r=2$.
Our next objective is to prove that every $k$-th order \dfo{} $\Tk$ in
$\tfS^{r\da}_{p,q}$, the subspace of $\tfS^r_{p,q}$ of \dfo s with
nonpositive $y$-degree, may be written as the sum of
a $k$-th degree polynomial in the generators of $\fg^{24,r}_p$ and an
operator in $\ker\pi^r_{p,q}$.
We shall construct a suitable collection of monomials in the
generators~\eqref{g24} such that their projections in $\tfS^r_{p,q}$
along $\ker\pi^r_{p,q}$ form a basis
of $\tfS^{r\da}_{p,q}$. (Note that in the triangular case every
monomial of degree at most
$n$ has a non-trivial projection in $\tfT_n$ along $\ker\pi_n$).
We follow the approach of the previous section,
obtaining first a set of monomials of maximal length (in this case the
bidegree does not determine uniquely the monomials of maximal length),
and then extend it to a basis of the image of the map
$\rho^{r,k}_p:\zS^k(\fgl_2\ltimes\CC^{r+1})\ra\fD^{(k)}$
determined by the representation
$\rho^r_p:\fgl_2\ltimes\CC^{r+1}\ra\fg^{24,r}_p$.
\begin{lemma}\label{lemma.bigbasis}\ \\[.3cm]
i) The monomials\footnote{By convention, the second group of monomials is not
present if $l=0$.}
$$
\big\{(J^{1+\ep}_p)^{l-j}(J^{5+\ep r})^j(J^3_p)^n(J^4)^{t-l-n},\,
(J^5)^{s}J^{5+i}(J^{5+r})^{l-s-1}(J^3_p)^n(J^4)^{t-l-n}\big\},
$$
with $\ep=0,1$, $0\leq j\leq l\leq t\leq k$, $0\leq s\leq l-1$,
$0\leq n\leq t-l$, and $1\leq i\leq r-\de_{s0}$ form a basis of
$\im\rho^{r,k}_p$ for all $k$.\\

\ni ii) The number of $k$-th order monomials in the basis is
$$
n_k=\frac 16(k+1)(k+2)\big((r+2)k+3\big)\,.
$$
\end{lemma}
The basis described in the previous lemma does not completely fulfill our
requirements, for it contains monomials of order less than $p$
acting trivially in $\cS^r_{p,q}$. A subset of
monomials acting non-trivially can be constructed by a simple graphical
argument, \cf Fig.~2.
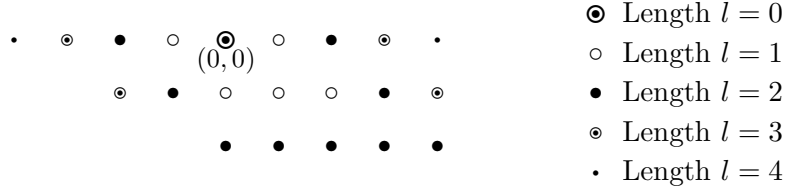
\begin{figure}[h]
\begin{center}
\begin{picture}(240,40)(-100,-40)
\put(0,0){\circle*{3}}\thicklines
\put(0,0){\circle{6}}\thinlines
\put(20,0){\circle{4}}\put(-20,0){\circle{4}}\put(0,-20){\circle{4}}
\put(20,-20){\circle{4}}\put(40,-20){\circle{4}}
\put(40,0){\circle*{4}}\put(-40,0){\circle*{4}}\put(-20,-20){\circle*{4}}
\put(0,-40){\circle*{4}}\put(20,-40){\circle*{4}}\put(40,-40){\circle*{4}}
\put(60,-40){\circle*{4}}\put(80,-40){\circle*{4}}\put(60,-20){\circle*{4}}
\put(60,0){\circle{4}}\put(-60,0){\circle{4}}
\put(60,0){\circle*{2}}\put(-60,0){\circle*{2}}
\put(-40,-20){\circle{4}}\put(80,-20){\circle{4}}
\put(-40,-20){\circle*{2}}\put(80,-20){\circle*{2}}
\put(80,0){\circle*{2}}\put(-80,0){\circle*{2}}
\put(140,10){\circle*{3}}\thicklines\put(140,10){\circle{6}}\thinlines
\put(150,7){\small Length $l=0$}
\put(140,-5){\circle{4}}\put(150,-8){\small Length $l=1$}
\put(140,-20){\circle*{4}}\put(150,-23){\small Length $l=2$}
\put(140,-35){\circle{4}}\put(140,-35){\circle*{2}}
\put(150,-38){\small Length $l=3$}
\put(140,-50){\circle*{2}}\put(150,-53){\small Length $l=4$}
\put(-11,-10){\footnotesize $(0,0)$}
\end{picture}
\end{center}
\caption{Lattice of bidegrees of monomials acting non-trivially
in $\cS^2_{4,2}$}
\end{figure}
\begin{lemma}\label{lemma.effbasis}
The projections in $\tfS^r_{p,q}$ along $\ker\pi^r_{p,q}$ of the monomials
\begin{equation}\label{effbasis}
\big\{(J^{1+\ep}_p)^{s}(J^{5+\ep r})^j(J^3_p)^n(J^4)^m,\,
(J^5)^{t}J^{5+i}(J^{5+r})^{j-t-1}(J^3_p)^n(J^4)^m\big\},
\end{equation}
with $\ep=0,1$, $0\leq j\leq q$, $0\leq s\leq p-jr$, $0\leq t\leq j-1$,
$1\leq i\leq r-\de_{t0}$, $0\leq n$, $0\leq m\leq q-j$, and
$n+r m\leq p-\tilde s-jr$ (where $\tilde s=s$ for the first group of monomials,
and $\tilde s=0$ for the second group), form a basis of $\tfS^{r\da}_{p,q}$.
\end{lemma}
{\em Proof.} The projections of the monomials in $\tfS^r_{p,q}$ along
$\ker\pi^r_{p,q}$ are certainly linearly independent. The number of monomials
in~\eqref{effbasis} is given by
$$
N^r_{pq}=\sum_{j=0}^q\sum_{s=0}^{p-jr}
f(j,s)\si\big(r,p-s-jr,\min\big(q,\big[\frac{p-s}r\big]\big)-j\big)\,,
$$
where $f(j,s)=2$ if $s\geq 1$ and $f(j,0)=jr+1$, and
$\si(r_0,p_0,q_0)=\dim\cS^{r_0}_{p_0,q_0}=(q_0+1)(p_0+1-\frac 12r_0q_0)$. We
leave as an exercise to the reader to verify that
$$
N^r_{pq}=\dim\tfS^{r\da}_{p,q}=
\frac{(q+1)(q+2)}2\Big((p+1)(p+1-qr)+\frac{qr^2}{12}(3q+1)\Big)\,.\QED
$$
The structure theorem for $\tfS^{r\da}_{p,q}$ now follows
from the fact that the highest
order derivatives of the monomials~\eqref{effbasis} are all independent.
\begin{thm}\label{thm.S}
Let $\Tk$ be a $k$-th order \dfo{} in $\tfS^{r\da}_{p,q}$. If $k\leq q$,
then $\Tk$ may be represented as
a $k$-th degree polynomial in the generators of $\fg^{24,r}_p$. If $k>q$,
then $\Tk$ may be expressed as the projection of such a polynomial in
$\tfS^r_{p,q}$ along $\ker\pi^r_{p,q}$.
\end{thm}

We shall now restrict ourselves to the ``untruncated'' modules
$\cS^r_p=\cS^r_{p,[p/r]}$ and the distinguished spaces of \dfo s
$\tfS^r_p=\tfS^r_{p,[p/r]}$ which preserve them.
We are interested in obtaining a method for constructing the most general
\dfo{} in $\tfS^r_p$ irrespective of their $y$-th degree. 
Following~\cite{Tu92b}, we introduce a well-adapted grading for the \dfo s in
$\tfS^r_p$. If $T\in\tfS^r_p$ has a well-defined bidegree $(i,j)$, we
define its {\em total degree} as $\Deg T=i+rj$.
The constructive method outlined in the proof of Lemma~\ref{lemma.iso}
allows us to obtain the \dfo{} $T$ representing any given endomorphism
of $\cS^r_p$. However, it does not give any information on the order of
$T$, which limits its value in quantum mechanical applications. The
following elementary lemma partially describes the subspace
\mbox{$\tfS^{r,(k)}_p=\tfS^r_p\cap\fD^{(k)}$} of $k$-th order
\dfo s.
\begin{lemma}
Let $T^{(k)}$ be an element of $\tfS^{r,(k)}_p$ homogeneous with
respect to the total degree. Then $-k r\leq\Deg T^{(k)}\leq k$.
\end{lemma}
{\em Proof.} The lower bound is obvious, for $T^{(k)}$ is a linear \dfo{} with
polynomial coefficients. Now, if $\Deg T^{(k)}=d>k$, $T^{(k)}$ must annihilate
the last $d$ diagonals $\{x^iy^j\:|\:i+rj=p-d+1,\ldots,p\}$, but a (nonzero)
$k$-th order~\dfo{} cannot annihilate more than $k$ diagonals.\QED\\

We can use this result and Theorem~\ref{thm.S} to construct the most
general \dfo{} in $\tfS^{r,(k)}_p$.
As an application, we give the explicit form of the most general
\dfo{} $T^{(2)}$ in $\tfS^{2,(2)}_p$,
which may be used to construct partially solvable \sch{} operators in
two variables by solving the corresponding equivalence problem, \cite{GKO94}.
We assume that $p\geq 4$ (otherwise some of the terms
would just have a non-vanishing projection on $\ker\pi^2_{p}$
along $\tfS^2_p$). Then, $T^{(2)}=\sum_{i=-4}^2\,T_i$, where:
\begin{align*}
& T_{-4} = a_1\,\Dy^2\,,\\
& T_{-3} = a_2\,\Dx\Dy+a_3\,x\Dy^2\,,\\
& T_{-2} = a_4\,\Dx^2+a_5\,x\Dx\Dy+a_6\,x^2\Dy^2+a_7\,y\Dy^2+a_8\,\Dy\,,\\
& T_{-1} = 
a_9\,x\Dx^2+a_{10}\,x^2\Dx\Dy+a_{11}\,y\Dx\Dy+a_{12}\,x^3\Dy^2+a_{13}\,xy\Dy^2
+a_{14}\,\Dx+a_{15}\,x\Dy\,,\\
& T_0 =
a_{16}\,x^2\Dx^2+a_{17}\,y\Dx^2+a_{18}\,x^3\Dx\Dy
+a_{19}\,xy\Dx\Dy+a_{20}\,x^4\Dy^2+a_{21}\,x^2y\Dy^2\\
&\qquad +a_{22}\,y^2\Dy^2+a_{23}\,x\Dx+a_{24}\,x^2\Dy+a_{25}\,y\Dy+a_{26}\,,\\
& T_1 = \big(a_{27}\,x^2\Dx+a_{28}\,y\Dx+a_{29}\,x
^3\Dy+a_{30}\,xy\Dy+a_{31}\,x\big)
        (x\Dx+2y\Dy-p)\,,\\
& T_2 = \big(a_{32}\,x^2+a_{33}\,y\big)(x\Dx+2y\Dy-p+1)(x\Dx+2y\Dy-p)\,.
\end{align*}
Note that $T_{-4},\ldots,T_0$ are linear combinations of monomials,
reflecting the fact that their total degree is nonpositive (an operator of
this type is sometimes referred to in the literature as
{\em exactly solvable} since it preserves an infinite flag of polynomial
subspaces). Note that the terms with coefficients $a_{17}$, $a_{28}$
and $a_{33}$ have positive $y$-degree, so they {\em cannot} be obtained from
a polynomial in the generators of $\fg^{24,2}_p$.
This is one of the first examples of a partially solvable \dfo{} preserving
a \fdim{} space of polynomials which is not Lie algebraic. (A number of
partially solvable operators which preserve
a \fdim{} module of monomials but are not Lie algebraic are also discussed
in~\cite{PT95}. The exponents of the monomials spanning these modules do
not form a continuous chain $\{0,1,\ldots,n\}$).

\section{Matrix $\mathbf k$-th order \dfo s in \boldmath$\CC^2$ with
invariant modules of polynomials}

We shall now consider the more general case of $k$-th order $2\times 2$
matrix \dfo s admitting \fdim{} invariant subspaces. Our goal is to prove
structure theorems for these operators, which are the counterparts of
Theorems~\ref{thm.R}, \ref{thm.T} and~\ref{thm.S},
in the case in which the invariant rank 2 module
$\bcN\subset\C^\infty(U,\CC\oplus\CC)$
is either a direct sum of rectangular, triangular or staircase modules.

Consider the associative algebra $\bfD=\bfD(U)$ of all $2\times 2$
matrix \dfo s with smooth coefficient functions. We introduce in $\bfD$
the usual $\ZZ_2$-grading; an operator
\begin{equation}\label{matrixT}
T=\bpm T_{11} & T_{12} \\ T_{21} & T_{22} \epm
\end{equation}
is said to be {\em even} if $T_{12}=T_{21}=0$, and {\em odd} if
$T_{11}=T_{22}=0$.
The associative algebra $\bfD$ becomes then a Lie superalgebra with
generalized Lie product given by
\begin{equation}\label{scomm}
[T_1,T_2]_s=T_1T_2-(-1)^{\deg T_1\deg T_2}T_2T_1\,.
\end{equation}
(We shall use the symbol $\bfD$ to denote both the associative algebra and
the Lie superalgebra). We will construct graded subalgebras of $\bfD$
of matrix \dfo s preserving the above mentioned vector-valued modules,
in analogy with the scalar case.

\subsection{The rectangular module $\cR_{n_1,m_1}\oplus\cR_{n_2,m_2}$}

This case has been recently studied in~\cite{BGK96}. We shall be very concise,
for the results are a natural generalization of the single-variable case
considered in~\cite{BK94} and~\cite{FGR96}. 
We shall assume that $n=n_2\geq n_1$ and $m=m_1\geq m_2$ (the other choice
$m_2\geq m_1$ leads to analogous results). Let $\De=n-n_1$ and $\Ga=m-m_2$.
Let us denote the direct sum $\cR_{n-\De,m}\oplus\cR_{n,m-\Ga}$ by
$\bcR^{n-\De,m}_{n,m-\Ga}$. Consider the graded subalgebra
$\fs^{\De,\Ga}_{n,m}$ of $\bfD$ generated by
the $8+2(\De+1)(\Ga+1)$ matrix \dfo s given by
\begin{alignat}{2}
& \Se=\diag(\Je{n-\De},\Je n)\,,& \qquad
& \Te=\diag(\Ke m,\Ke{m-\Ga})\,,\qquad \ep=\pm,0\,,\notag\\
& J=\frac 12\diag(n+\De,n)\,,& \qquad & K=\frac 12\diag(m,m+\Ga)\,,\notag\\
& \Qab^-=q_\al(x)\,\bq_\be(y,m,\Ga)\,\si^-\,,
& \qquad & \Qab^+=\bq_\al(x,n,\De)\,q_\be(y)\,\si^+\,,\label{sDeGa}
\end{alignat}
where $\al=0,\ldots,\De$\,,\; $\be=0,\ldots,\Ga$\,,\;
$\si^+=(\si^-)^t=\bpm 0 & 1 \\ 0 & 0\epm$, and
$$
q_\al(x)=x^\al\,,\qquad \bq_\al(x,n,\De)=
\prod_{k=1}^{\De-\al}(x\Dx-n+\De-k)\,\Dx^\al\,,
$$
with the convention that a product with its lower limit greater than
the upper one is defined to be $1$. Let $\fs^{\De,\Ga}$ be the abstract
Lie superalgebra corresponding to the realization $\fs^{\De,\Ga}_{n,m}$
(it may be shown that the supercommutators~\eqref{scomm}
of the operators~\eqref{sDeGa} do not depend explicitly on
$n$ and $m$, \cite{BK94}). Just as in the
single-variable case, \cite{FGR96}, the Lie superalgebra $\fs^{\De,\Ga}$
will be infinite-dimensional unless the gaps $\De$ and $\Ga$ are suitably
restricted. The corresponding \fdim{} Lie superalgebras are easily identified.
\begin{lemma}
The Lie superalgebra $\fs^{\De,\Ga}$ is \fdim{} if and only if
$\De+\Ga\leq 1$. We have:
$$
\fs^{1,0}\iso\fs^{0,1}\iso\fspl(2,1)\oplus\fsl_2\,,\qquad
\fs^{0,0}\iso\overset{3}{\underset{i=1}\oplus}\,\fsl^i_2
$$
\end{lemma}
{\em Remarks.} We have not included in $\fs^{1,0}$ the central term given
by $K$. We regard $\fs^{0,0}$ as a Lie algebra with $J$ and $K$ replaced by
$\tJ=\si_3$. Finally, recall the isomorphism $\fspl(2,1)\iso\fosp(2,2)$,
\cite{Sc79}.\\

Let $\bfR^{n-\De,m}_{n,m-\Ga}$ be the associative subalgebra of $\bfD$
spanned by the operators which preserve $\bcR^{n-\De,m}_{n,m-\Ga}$.
It is easy to see that $\fs^{\De,\Ga}_{n,m}\subset\bfR^{n-\De,m}_{n,m-\Ga}$,
and that the action is irreducible. Thus, the Burnside Theorem guarantees
that any matrix \dfo{} preserving this module may
be expressed as a polynomial in the generators of $\fs^{\De,\Ga}_{n,m}$
plus an operator annihilating the module. The kernel of the homomorphism
$\pi^{\De,\Ga}_{n,m}:\bfR^{n-\De,m}_{n,m-\Ga}\ra\End(\bcR^{n-\De,m}_{n,m-\Ga})$
is simply the set of operators~\eqref{matrixT} whose component
$T_{ij}$ belongs to $\ker\pi_{n_j,m_j}$ as given in Lemma~\ref{lemma.ker}.
We choose a distinguished complement $\tbfR^{n-\De,m}_{n,m-\Ga}$ to
$\ker\pi^{\De,\Ga}_{n,m}$ just as in the scalar case.
The following lemma is a direct generalization of the corresponding
single-variable result, \cite{FGR96}:
\begin{lemma}
Let $\brho^{\De,\Ga}_{n,m}:\fU(\fs^{\De,\Ga})\ra\bfD$ be the homomorphism
determined by the representation
$\rho^{\De,\Ga}_{n,m}:\fs^{\De,\Ga}\ra\fs^{\De,\Ga}_{n,m}$. The monomials
\begin{multline}\label{basisR}
\big\{X(S^\pm)^i(S^0)^s(T^\pm)^j(T^0)^t\,,\;
Q^\vep_{\al 0}(S^{\vep^\ast})^i(T^\pm)^j(T^0)^t\,,\\
Q^\vep_{0 \be}(S^\pm)^i(S^0)^s(T^\vep)^j\,,\;
\Qab^\vep(S^{\vep^\ast})^i(T^\vep)^j\big\}\,,
\end{multline}
where $X=Q^\vep_{00}$, $\id$, $\tJ$, $\vep=-\vep^\ast=\pm$,
$\al=1,\ldots,\De$, and $\be=1,\ldots,\Ga$,
form a basis of $\im\brho^{\De,\Ga}_{n,m}$.
\end{lemma}
Note that $\tJ$ may be expressed in terms of $\id$, $J$, $K$
(or $Q^\pm_{00}$ for $\fs^{0,0}$).
The number of monomials in~\eqref{basisR} acting non-trivially
in $\bcR^{n-\De,m}_{n,m-\Ga}$ is precisely
$\dim\tbfR^{n-\De,m}_{n,m-\Ga}$.
\begin{thm}
Let $\Tk$ be a $k$-th order \dfo{} in $\tbfR^{n-\De,m}_{n,m-\Ga}$.
Then $\Tk$ may be expressed as a linear combination of the monomials
in~\eqref{basisR} of differential order at most $k$.
\end{thm}

\subsection{The triangular module $\cT_{n_1}\oplus\cT_{n_2}$}

We shall assume that $n=n_2\geq n_1$, and let $\De=n-n_1$. Let
us denote the module $\cT_{n-\De}\oplus\cT_n$ by $\bcT^{n-\De}_n$.
Let $\bfT^{n-\De}_n$ be the associative subalgebra of
$\bfD$ of \dfo s preserving $\bcT^{n-\De}_n$. Consider the
$9+(\De+1)(\De+2)$ matrix \dfo s given by
\begin{align}
& T^i=\diag(J^i_{n-\De},J^i_n)\,,\qquad J=
\frac 13\diag(n+2\De,n)\,,\qquad i=1,\ldots,8\,,\notag\\
& \Qab^-=q_{\al\be}\,\si^-\,,\qquad \Qab^+=\bq_{\al\be}(n,\De)\,\si^+\,,\qquad
0\leq\al+\be\leq\De\,,\label{sDe}
\end{align}
with $J^i_n$ given in~\eqref{g15}, and
$$
q_{\al\be}=x^\al y^\be\,,\qquad
\bq_{\al\be}(n,\De)=
\prod_{k=1}^{\De-\al-\be}(x\Dx+y\Dy-n+\De-k)\,\Dx^\al\,\Dy^\be\,.
$$
We observe that the even generators in~\eqref{sDe} span a Lie algebra
which is isomorphic to $\fgl_3$. It may be easily verified that both
the generators $\Qab^-$ and $\Qab^+$
span irreducible modules under the adjoint action of the even subalgebra,
which are isomorphic to the standard cyclic $\fsl_3$-module of highest
weight $(\De,0)$, with $J$ acting as a scalar. Moreover, the
anticommutator of $\Qab^-$ and $Q^+_{\ga\de}$
can be expressed as a polynomial of degree $\De$ in the even generators and 
the identity, depending explicitly on $n$ only through a multiple of the
identity. Let us denote by $\fs^\De_n$ the Lie superalgebra
generated by the \dfo s~\eqref{sDe}.
\begin{lemma}
The Lie superalgebra $\fs^\De_n$ is \fdim{} if and only if $\De\leq 1$. We have:
$$
\fs^1_n\iso\fspl(3,1)\,,\qquad\fs^0_n\iso\fsl_3\oplus\fsl_2\,.
$$
\end{lemma}
{\em Remark.} We regard $\fs^0_n$ as a Lie algebra with $J$ replaced by
$\tJ=\si_3$.

In general, the Lie superalgebra of \dfo s preserving a pair of
$N$-dimensional ``triangular'' modules with a gap $\De=1$ is isomorphic to
$\fspl(N+1,1)$. It is easy to see that $\fs^\De_n\subset\bfT^{n-\De}_n$,
and that the action on $\bcT^{n-\De}_n$ is irreducible. The kernel of
the homomorphism $\pi^\De_n:\bfT^{n-\De}_n\ra\End(\cT^{n-\De}_n)$
is just the set of the operators~\eqref{matrixT} with entries $T_{ij}$ 
in $\ker\pi_{n_j}$ as given in Lemma~\ref{lemma.ker}. As usual,
we choose as distinguished complement to $\ker\pi^\De_n$ the
vector space $\tbfT^{n-\De}_n$ spanned by the \dfo s in $\bfT^{n-\De}_n$
with components of the form
$T_{ij}=\sum_{s+t\leq n_j} f_{st}(x,y)\,\Dx^s\Dy^t$.\\

We first note that the definition of bidegree can be naturally extended to
matrix-valued \dfo s
acting on $\CC[x,y]\oplus\CC[x,y]$. The bidegrees of the generators~\eqref{sDe}
of $\fs^\De_n$ are well-defined. The bidegree of each $T^i$ is given by
the bidegree of the corresponding $J^i_n$, and $\deg J=(0,0)$. For the
generators of the odd subspace,
we have $\deg Q^\pm_{\al\be}=\mp(\al,\be)$. In order to define the length
of a monomial in the generators~\eqref{sDe} we need to be a little more
careful. Let $T^K=(T^1)^{k_1}\cdots (T^8)^{k_8}$.
Consider the monomials of the form $T^K_\ep=X_\ep\,T^K$, where $\ep=1,2,\pm$, 
$X_1=P_1$, $X_2=P_2$, $X_\pm=Q^\pm_{\al\be}$, and $P_1$, $P_2$ are the
canonical projectors, which can of course be expressed in terms of
$\id$ and $J$ (or $\tJ$, if $\De=0$). We define the length of
a monomial $T^K_\ep$ as
$$
|T^K_\ep|=\min_{\deg T^L_\ep=\deg T^K_\ep}|L|\,,
$$
where the minimum is taken over the set of monomials $T^L_\ep$ with
$|L|\leq |K|$. Again, a monomial $T^K_\ep$ is of maximal length if
$|T^K_\ep|=|K|$. Thus, $P_1$, $P_2$ and the $Q^\pm_{\al\be}$'s are
all of maximal length (equal to zero).
For $\ep=\pm$, the bidegree does {\em not} determine uniquely the
monomials of maximal length.
\begin{figure}[h]
\begin{center}
\begin{picture}(240,100)(-120,-60)
\put(0,0){\circle*{3}}\thicklines\put(0,0){\circle{6}}\thinlines
\put(-20,0){\circle*{3}}\thicklines\put(-20,0){\circle{6}}\thinlines
\put(0,-20){\circle*{3}}\thicklines\put(0,-20){\circle{6}}\thinlines
\put(20,0){\circle{4}}\put(0,20){\circle{4}}\put(-20,20){\circle{4}}
\put(-40,20){\circle{4}}\put(-40,0){\circle{4}}\put(-20,-20){\circle{4}}
\put(0,-40){\circle{4}}\put(20,-40){\circle{4}}\put(20,-20){\circle{4}}
\put(40,0){\circle*{4}}\put(20,20){\circle*{4}}\put(0,40){\circle*{4}}
\put(-20,40){\circle*{4}}\put(-40,40){\circle*{4}}\put(-60,40){\circle*{4}}
\put(-60,20){\circle*{4}}\put(-60,0){\circle*{4}}\put(-40,-20){\circle*{4}}
\put(-20,-40){\circle*{4}}\put(0,-60){\circle*{4}}\put(20,-60){\circle*{4}}
\put(40,-60){\circle*{4}}\put(40,-40){\circle*{4}}\put(40,-20){\circle*{4}}
\put(120,20){\circle*{3}}\thicklines\put(120,20){\circle{6}}\thinlines
\put(130,17){\small Length $l=0$}
\put(120,0){\circle{4}}\put(130,-3){\small Length $l=1$}
\put(120,-20){\circle*{4}}\put(130,-23){\small Length $l=2$}
\put(-11,-10){\footnotesize $(0,0)$}
\end{picture}
\end{center}
\caption{Lattice of bidegrees of monomials $T^K_+$ of length $l\leq 2$ for
$\De=1$.}
\end{figure}
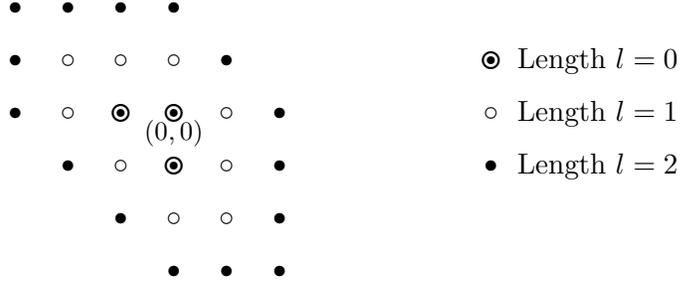
The structure theorem for $\tbfT^{n-\De}_n$ may now be stated as follows:
\begin{thm}
Let $\Tk$ be a $k$-th order \dfo{} in $\tbfT^{n-\De}_n$. Then $\Tk$
may be expressed as a linear combination of the monomials
\begin{equation}\label{sDebasis}
\big\{T^L_\ep\,(T^7)^i(T^8)^j\;\big|\;
\ep=1,2,\pm\,,\; 0\leq i+j+|L|\leq t\leq k-\de_{\ep+}\De\big\}\,,
\end{equation}
where $T^L_\ep$ are monomials of maximal length and
different bidegrees for each (fixed) value of $\ep$.
\end{thm}
{\em Remark.} If $\De>k$, no monomial $T^L_+\,(T^7)^i(T^8)^j$ may appear in
the expansion of $\Tk$.

{\ni\em Proof.} The monomials~\eqref{sDebasis} are linearly independent.
Let $N_k^\ep$ denote the number of monomials of order at most $k$ and type
$\ep$ in~\eqref{sDebasis}. We already know that
$N_k^{1,2}=\frac 14(k+1)^2(k+2)^2$, \cf \eqref{dimim}. Let
$\tau(n_0)=\dim\cT_{n_0}=\frac 12(n_0+1)(n_0+2)$. An easy graphical argument,
\cf Fig.~3, shows that
\begin{align}
N_k^\pm &= \frac{(\De+1)(\De+2)}2\,\tau(k-\de_{\ep+}\De)
            +\sum_{l=1}^{k-\de_{\ep+}\De}3(2l+\De)\,
            \tau(k-\de_{\ep+}\De-l)\notag\\
        &= \frac 14(k+1)(k+2)(k\mp\De+1)(k\mp\De+2)\,.\label{nkpm}
\end{align}
Therefore, $N_{n-\De}^1+N_n^2+N_n^++N_{n-\De}^-=\dim\tbfT^{n-\De}_n$.
The theorem follows from the fact that the highest order derivatives
of the monomials~\eqref{sDebasis} are all independent.\QED\\

The dimension of the subspace $\tbfT^{n-\De,(k)}_n$ of $\tbfT^{n-\De}_n$
of \dfo s of order at most $k$ is given by
$$
\dim\tbfT^{n-\De,(k)}_n=\frac 12(k+1)(k+2)(2k^2+6k+4+\De^2)\,.
$$
It is worth mentioning that in the one-dimensional rank two case this
number is independent of $\De$, \cite{FGR96}.

\subsection{The staircase module $\cS^r_{p_1,q_1}\oplus\cS^r_{p_2,q_2}$}

For the sake of simplicity, we shall restrict ourselves to the
case $p=p_1+1=p_2$, $q=q_1=q_2$. Let us denote the module
$\cS^r_{p-1,q}\oplus\cS^r_{p,q}$ by $\bcS^{r,1}_{p,q}$.
Consider the \fdim{} Lie superalgebra $\fs^r_p$ spanned by the $10+2r$
matrix \dfo s given by:
\begin{align}
& T^i=\diag(J^i_{p-1},J^i_p)\,,\qquad J=\frac 12\diag(p+1,p)
  \,,\qquad i=1,\ldots,5+r\,,\notag\\
& \Qm_0=\si^-\,,\qquad \Qm_1=x\,\si^-\,,\qquad
  \Qp_0=(x\Dx+ry\Dy-p)\,\si^+\,,\qquad
  \Qp_1=\Dx\,\si^+\,,\notag\\[2mm]
& \Qp_{2+\al}=x^\al\Dy\,\si^+\,,\qquad \al=0,\ldots,r-1\,,\label{sr}
\end{align}
with $J^i_p$ given in~\eqref{g24}. The abstract structure of the Lie
superalgebras $\fs^r_p$ is the semidirect product
$\fpl(2,1)\ltimes\CC^{r+1,r}$, where $\CC^{r+1,r}$ corresponds
to the Abelian Lie superalgebra spanned by $T^{5+i}$ and $\Qp_{2+\al}$, 
$i=0,\ldots,r$, $\:\al=0,\ldots,r-1$. It may be readily verified that
$\fs^r_p$ is contained in the associative subalgebra $\bfS^{r,1}_{p,q}$
of $\bfD$ of \dfo s preserving $\bcS^{r,1}_{p,q}$. 
The kernel of the homomorphism
$\pi^{r,1}_{p,q}:\bfS^{r,1}_{p,q}\ra\End(\bcS^{r,1}_{p,q})$
is spanned by the \dfo s~\eqref{matrixT} with components $T_{ij}$ in
$\ker\pi^r_{p_j,q}$ as given in Lemma~\ref{lemma.ker}. We choose as
a distinguished complement to $\ker\pi^{r,1}_{p,q}$
the subspace $\tbfS^{r,1}_{p,q}\subset\bfS^{r,1}_{p,q}$ of \dfo s with
components of the form:
$$
T_{ij}=\sum_{\substack{0\leq s+rt\leq p_j\\0\leq t\leq q}}
f_{st}(x,y)\,\Dx^s\Dy^t\,.
$$
We note that just as in the scalar case, the module $\bcS^{r,1}_{p,q}$
is reducible (but not completely reducible) under the action of $\fs^r_p$.
In fact, no operator in $\tbfS^{r,1}_{p,q}$ with positive $y$-degree
may be obtained from a polynomial in the generators of $\fs^r_p$. Let
$\tbfS^{r,1\da}_{p,q}$ denote the subspace of $\tbfS^{r,1}_{p,q}$ of
\dfo s with nonpositive $y$-degree. The structure
theorem for $\tbfS^{r,1\da}_{p,q}$ may be proved by essentially the same
argument used in Theorem~\ref{thm.S}.
\begin{thm}
Let $\Tk$ be a $k$-th order \dfo{} in $\tbfS^{r,1\da}_{p,q}$.
If $k\leq q$, then $\Tk$ may be represented as
a $(k+1)$-th degree polynomial in the generators of $\fs^r_p$. If $k>q$,
then $\Tk$ may be expressed as the projection of such a polynomial in
$\tbfS^{r,1}_{p,q}$ along $\ker\pi^{r,1}_{p,q}$.
\end{thm}
{\em Proof.} The key of the proof is to obtain a suitable set of
monomials in the generators of $\fs^r_p$ such that their projections in
$\tbfS^{r,1}_{p,q}$ along $\ker\pi^{r,1}_{p,q}$ are all independent.
We consider monomials of the form
\begin{equation}\label{basisS}
T^L_\ep\,(T^3)^i(T^4)^j,
\end{equation}
where $\ep=1,2\pm$, and $T^L_\ep$ are monomials of maximal
length and different bidegrees for any (fixed) value of $\ep$.
For $\ep=1,2$, Lemma~\ref{lemma.effbasis}
implies that we can obtain $\dim\tfS^{r\da}_{p_\ep,q}$ monomials of
the form~\eqref{basisS} with independent
projections in $\tbfS^{r,1}_{p,q}$ along $\ker\pi^{r,1}_{p,q}$. Let
$\tfS^{r\da-}_{p,q}$ denote the subset of $\tbfS^{r,1\da}_{p,q}$ of \dfo s
mapping $P_1(\bcS^{r,1}_{p,q})$ into $P_2(\bcS^{r,1}_{p,q})$.
The projections in $\tbfS^{r,1}_{p,q}$ along $\ker\pi^{r,1}_{p,q}$
of the monomials
\begin{equation}\label{basisS-}
\big\{\Qm_\ep(T^{1+\ep})^s(T^{5+\ep r})^j(T^3)^n(T^4)^m,\,
\Qm_0(T^5)^{t}T^{5+i}(T^{5+r})^{j-t-1}(T^3)^n(T^4)^m\big\}\,,
\end{equation}
with $\ep=0,1$, $0\leq j\leq q$, $0\leq s\leq p-1-jr$, $0\leq t\leq j-1$,
$1\leq i\leq r-\de_{t0}$, $0\leq n$, $0\leq m\leq q-j$, and
$n+r m\leq p-\tilde s-jr$ (where $\tilde s=s$ for the first group
of monomials, and $\tilde s=0$
for the second group), form a basis of $\tfS^{r\da-}_{p,q}$.
The number of monomials in~\eqref{basisS-} is given by
$$
N^{r-}_{pq}=\sum_{j=0}^q\sum_{s=0}^{p-1-jr}
f(j,s)\si\big(r,p-1-s-jr,\min\big(q,\big[\frac{p-1-s}r\big]\big)-j\big)\,,
$$
where $f(j,s)=2$ if $s\geq 1$ and $f(j,0)=jr+2$, and
$\si(r_0,p_0,q_0)=\dim\cS^{r_0}_{p_0,q_0}$. It may be verified that
$$
N^{r-}_{pq}=\dim\tfS^{r\da-}_{p,q}=
\frac{(q+1)(q+2)}2\Big(p(p+1)-qr(p+\frac 23)+\frac{qr^2}{12}(3q+1)\Big)\,.
$$
We also leave as an exercise to the reader to obtain
$$
N^{r+}_{pq}=\dim\tfS^{r\da+}_{p,q}=
\frac{(q+1)(q+2)}2\Big(p(p+1)-qr(p+\frac 13)+\frac{qr^2}{12}(3q+1)\Big)
$$
monomials of the form $T^L_+\,(T^3)^i(T^4)^j$ with independent projections in
$\tbfS^{r,1}_{p,q}$ along $\ker\pi^{r,1}_{p,q}$. The theorem follows from the
fact that the highest order derivatives
of the monomials~\eqref{basisS} are all different.\QED

\vspace{.5cm}

\section*{Acknowledgements}

It is a pleasure to thank Artemio Gonz\'alez-L\'opez and
Miguel A. Rodr\'\i guez for helpful discussions. One of us (F.F.)
wishes to thank the Fields Institute for its hospitality.

\end{document}